\documentclass{article}
\usepackage[latin1]{inputenc}
\usepackage{graphicx}
\usepackage{amsmath, enumerate, amsfonts}
\usepackage[widemargins]{a4}
\usepackage[all]{xy}
\xyoption{arc}

\newtheorem{Theorem}{Theorem}[section]

\title{Loop homology of spheres and complex projective spaces}
\author{Nora Seeliger}

\begin{document}
\maketitle
\begin{abstract}
In his Inventiones paper Ziller (1977) computed the integral homology as a graded abelian group of the free loop space of compact, globally symmetric spaces of rank 1. Chas and Sullivan (1999)showed that the homology of the free loop space of a compact closed orientable manifold can be equipped with a loop product and a BV-operator making it a Batalin-Vilkovisky algebra.
Cohen, Jones and Yan developed (2004) a spectral sequence
which converges to the loop homology as a spectral sequence of algebras. They computed the algebra structure of the loop homology of spheres and complex projective spaces by using Ziller's results and the method of Brown-Shih (1959, 1962). In this note we compute the loop homology algebra by using only spectral sequences and the technique of universal examples. We therefore not only obtain Ziller's and Brown-Shih's results in an elementary way, but we also replace the roundabout computations of Cohen, Jones and Yan (2004) making them independent of Ziller's and Brown-Shih's work. Moreover we offer an elementary technique which we expect can easily be generalized and applied to a wider family of spaces, not only the globally symmetric ones.
\end{abstract}
\section{Introduction}
Ziller studied in \cite{[10]} the free loop space of globally symmetric spaces of rank 1 which - due to the classification of Cartan, see \cite{Cartan1} and \cite{Cartan2} following the discussion in \cite{Helgason} - are spheres, projective spaces and the Cayley projective plane. Using the Energy function as a Morse function he computed the integral homology as a graded abelian group with the help of Morse theory,  methods from differential geometry and by rigorously using the symmetries of the spaces under consideration. The results of Brown and Shih, see \cite{Brown} and \cite{Shih} or the discussion in \cite{McCleary} allow to compute the differentials in the Serre spectral sequence for the free loop space fibration in the case that the base space is $(n-1)-$connected via the transgression homomorphism and the Hopf algera structure of the homology of the fibre. In \cite{[8]} Smith computed the characteristic $0$ cohomology of the free loop space for a complex projective space by making use of the Eilenberg Moore spectral sequence, the Eilenberg Moore theorem and constructing an explicit Koszul-type resolution. In \cite{[7]} we only made use of an elementary computation with the Serre spectral sequence. A related result of Smith in charasteristic $p$ appeared in \cite{LS2}. In \cite{[7]+} we described how Fadell and Husseini in \cite{FadellHusseini} computed the cohomology of the free loop space of a complex projective which is discussed as well in \cite{CrabbJames}. For the sake of completeness we mention that Westerland \cite{[9]} computed the string homology of spheres and complex projective spaces in \cite{[9]}, and Menichi and Hepworth computed the Batalin-Vilkovisky algebra structure for spheres respectively complex projective spaces in \cite{[5]} respectively \cite{[4]}. For simply connected spaces $X$ and commutative rings $k$ Goodwillie, Burghlea and Fiedorowicz proved that the Hochschild cohomology of the singular chains on the space of pointed loops $HH^*S_*(\Omega X)$ is isomorphic to the free loop spaces cohomology $H^*(\Lambda (X))$. Menichi proved in \cite{[6]} that this isomorphism is compatible with the usual cup product on $H^*(\Lambda (X))$ and the cup product of Cartan and Eilenberg on $HH^*S_*(\Omega X)$ and he used this isomorphism to compute the cohomology of the free loop space for a suspended space, for a complex projective space and a finite $CW-$space of dimension $p$ such that $(p-1)! \in k^{\times}$. A stable splitting of free loop spaces of spheres and complex projective spaces by Boekstedt and Ottosen has been discussed in \cite{BO}.\\ In \cite{[2]} Cohen, Jones and Yan computed the loop homology algebra of spheres and projective spaces by constructing a spectral sequence and then working backwards from Ziller's results for the case of complex projective spaces and used the method of Brown and Shih, see \cite{Brown} and \cite{Shih} for the differentials in the loop homology algebra spectral sequence for spheres which we will sketch as described in \cite{McCleary} as we go along.\\ In this note we compute the loop homology (as an algebra) independently of all these results, only using elementary spectral sequence computations, basic properties of the loop product of Chas and Sullivan, and the technique of universal examples by which we mean to map the fibration for which we want to compute the Serre spectral sequence to another fibration for which we understand the Serre spectral sequence better in order to use the naturality property to be able to compute differentials in the original spectral sequence. This paper therefore not only reproduces Ziller's results in an elementary way. It replaces as well the roundabout computation of Cohen, Jones and Yan and moreover we offer an elementary technique which we expect can easily be generalized and applied to a wider family of spaces, not only the globally symmetric ones.  \\
In this article $M$ always denotes a compact, closed, oriented manifold and let $\Lambda (M)=map(S^1,M)$ be its free loop space i. e.  %$\mathbb{C}P^n$ is the projective space of the vector space $\mathbb{C}^{n+1}$ of dimension $n+1$ over the complex numbers $\mathbb{C}$ and $\Lambda(\mathbb{C}P^n)$ is 
the function space $map(S^1,M)$ of unpointed maps from a circle $S^1$ into $M$ topologized with the compact open topology. Homology and cohomology will always be understood with integer coefficients. Following standard notation $\mathbb{H}_*(\Lambda (M))$ denotes the loop homology algebra in the sense of Chas and Sullivan for a compact oriented closed manifold $M$, the functor $E(-)$ denotes the exterior algebra functor whereas $\mathbb{Z}[-]$ is a polynomial algebra with integer coefficients and $\Gamma[\gamma]$ is a divided power algebra.\\
This work is partly based on the author's diploma thesis at the Georg-August-Universität Göttingen under the supervision of Prof. Dr. Thomas Schick. The author was supported by DFG grant $9209212$ "Smooth Cohomology". I would like to thank my advisor, Prof. Dr. Ulrich Bunke, and the referee.\\
We prove the following theorems.
\begin{Theorem} There is an isomorphism of algebras \begin{eqnarray*}\mathbb{H}_*(\Lambda (S^1);\mathbb{Z})\cong\mathbb{Z}[t, t^{-1}]\otimes E(x)\end{eqnarray*} for $\mathbb{H}_*(\Lambda (S^1);\mathbb{Z})$ the integral loop homology algebra of $S^1$ with $|x|=-1$, and $|t|=0$.
\end{Theorem}
\begin{Theorem}For $n>1$ odd there is an isomorphism of algebras \begin{eqnarray*}\mathbb{H}_*(\Lambda (S^n);\mathbb{Z})\cong\mathbb{Z}[y]\otimes E(x)\end{eqnarray*} for $\mathbb{H}_*(\Lambda (S^n);\mathbb{Z})$ the integral loop homology algebra of odd dimensional spheres with $|x|=-n$, and $|y|=n-1$.
\end{Theorem}
\begin{Theorem}There is an isomorphism of algebras 
\begin{eqnarray*}\mathbb{H}_*(\Lambda (\mathbb{C}P^n));\mathbb{Z}) \cong\frac{ E(w)\otimes
 \mathbb{Z}[x, y]}{(x^{n+1}, (n + 1)x^ny,w\cdot x^n)},\end{eqnarray*} for $\mathbb{H}_*(\Lambda (\mathbb{C}P^n));\mathbb{Z})$ the integral loop homology algebra of complex projective spaces with $|x|=-2$, $|y|=2n$, and $|w|=-1$.
\end{Theorem}
\begin{Theorem}For $n\geq 2$ even there is an isomorphism of algebras 
\begin{eqnarray*}\mathbb{H}_*(\Lambda (S^n));\mathbb{Z}) \cong\frac{E(z)\otimes\mathbb{Z}[x, y]}{( x^2, x\cdot z,  2x\cdot y)},
\end{eqnarray*} 
for the integral loop homology algebra $\mathbb{H}_*(\Lambda (S^n));\mathbb{Z})$ with 
\begin{eqnarray*}
|x|=-n\text{, }|y|=2n-2\text{, and }|z|=-1.
\end{eqnarray*} 
\end{Theorem}
%\textbf{Theorem 4:} \\
We begin with quoting a theorem from \cite{[3]} which characterizes the loop homology algebra.
\begin{Theorem}
Let $M$ be a compact, closed, oriented manifold. Then the loop product defines a map
\begin{eqnarray*}
\mu_*:H_*(\Lambda(M))\otimes H_*(\Lambda(M))\rightarrow H_*(\Lambda(M)),
\end{eqnarray*}
\textit{making $\mathbb{H}_*(\Lambda(M))$ an associative, graded commutative algebra.}
\end{Theorem}
Moreover we have from \cite{[2]} the following spectral sequence.
\begin{Theorem}
There is a second quadrant spectral sequence of algebras 
\begin{eqnarray*}
\{\mathbb{E}^r_{p,q},d^r|p\leq 0,q \geq 0\}
\end{eqnarray*}
such that
\begin{enumerate}
\item $\mathbb{E}^r_{*,*}$ is an algebra and the differential $d^r:\mathbb{E}^r_{*,*}\rightarrow \mathbb{E}^r_{*-r,*+r-1}$ fulfills the Leibniz rule for each $r\geq 1$.
\item The spectral sequence converges to $\mathbb{H}_*(\Lambda(M))$ as a spectral sequnce of algebras. 
\item For $m,n\geq 0$ we have
\begin{eqnarray*}
\mathbb{E}^2_{-m,n}\cong H^m(M,H_n(\Omega (M)).
\end{eqnarray*}
Further, the isomorphism $\mathbb{E}^2_{-*,*}\cong H^*(M;H_*(\Omega(M)))$ is an isomorphism of algebras, where the algebra structure of $H^*(M;H_*(\Omega(M)))$ is given by the cup product on the cohomology of $M$ with coefficients in the Pontryagin ring $H_*(\Omega(M))$.
\item The spectral sequence is natural with respect to smooth maps between manifolds.
\end{enumerate}
\end{Theorem}
We continue comparing the Chas-Sullivan product $(\mathbb{H}_*(M), \bullet)$ with the usual homology of the manifold with intersection product $(H_*(M), \wedge)$ and the homology of the based loop space with Pontryagin product $(H_*(\Omega(M)), \cdot)$ (see \cite[Section 3]{[1]}). We have two maps,
\begin{eqnarray*}
H_*(M)\overset{\epsilon}{\longrightarrow }\mathbb{H}_*(M)\overset{\cap}{\longrightarrow }H_*(\Omega(M))
 \end{eqnarray*}
where $\epsilon$ is the inclusion of constant loops into all loops and $\cap$ is the transversal intesection with one fiber of the projection loop space $\overset{evaluation}{\longrightarrow }M$. If we use the usual grading on $H_*(\Omega(M))$, the shifted grading on $\mathbb{H}_*(M)$, and the analogous shifted grading on $H_*(M)$, then these products and the two maps have degree zero. The composition
\begin{eqnarray*}
H_*(M)\overset{\epsilon}{\longrightarrow }\mathbb{H}_*(M)\overset{\cap}{\longrightarrow }H_*(\Omega(M))
 \end{eqnarray*}
preserves products as follows from the definitions. Moreover $\epsilon$ is an injection onto a direct summand, and for any Lie group manifold, $\cap$ is a surjection.\\
With this input we can start the computations by looking at the fibration $\Omega(M)\rightarrow\Lambda(M)\overset{eval}{\longrightarrow }M$ where $eval :\Lambda (M)\rightarrow M$ denotes the evaluation at the basepoint of $S^1$. The strategy of the proofs will be to compare the loop homology algebra spectral sequence $\{\mathbb{E}^r,d^r\}$ and the Serre spectral sequence.
\section{Odd dimensional spheres}
Before we treat the general case we give an argument which deals with the case $S^1$. Recall that $S^1$ is an H-space and therefore the free loop space $\Lambda (S^1)$ is as well. The following argument is taken from \cite{LS2}, Lemma 3. With the H-space structure discussed above the map $eval :\Lambda (S^1)\rightarrow S^1$ denoting the evaluation at the basepoint of $S^1$ becomes an $H$-map making \begin{eqnarray*}\Omega(S^1)\rightarrow\Lambda(S^1)\overset{eval}{\longrightarrow }S^1\end{eqnarray*} a principal bundle which has a cross section $s:S^1\rightarrow\Lambda (S^1), x\mapsto$ constant loop at $x$. Since a principal fibration with cross section is trivial we obtain a homotopy equivalence \begin{eqnarray*}\Lambda(S^1)\simeq \Omega(S^1)\times S^1\simeq\mathbb{Z}\times S^1\end{eqnarray*} since all the components of $\Omega(S^1)$ are contractible. The loop homology algebra structure is given by the intersection product structure on $H_*(S^1)$ with the group algebra structure \begin{eqnarray*}H_0(\mathbb{Z})\cong \mathbb{Z}[t,t^{-1}]\end{eqnarray*}. With Poincaré duality we obtain our first theorem where the class $x\in\mathbb{H}_{-1}(\Lambda (S^1);\mathbb{Z})$ corresponds to the generator in $H^1(S^1;\mathbb{Z})$. It follows from the discussion at the end of the introduction that this multiplicative structure is the Chas-Sullivan product structure because $H_*(S^1)$ injects onto a direct summand of $\mathbb{H}_*(\Lambda (S^1))$ and $\mathbb{H}_*(\Lambda (S^1))$ surjects onto $H_*(\Omega (S^1))$ since $S^1$ is a Lie group.
\begin{Theorem}There is an isomorphism of algebras \begin{eqnarray*}\mathbb{H}_*(\Lambda (S^1);\mathbb{Z})\cong\mathbb{Z}[t,t^{-1}]\otimes E(x)\end{eqnarray*} for $\mathbb{H}_*(\Lambda (S^1);\mathbb{Z})$ the integral loop homology algebra of $S^1$ with $|x|=-1$, and $|t|=0$.
\end{Theorem}
For the loop homology of $S^n$ where $n>1$ is odd consider the fibration 
$\Omega(S^n)\hookrightarrow\Lambda(S^n)\stackrel{eval}{\rightarrow}S^n$
and recall $H^*(S^n)\cong E(x)$, $H_*(\Omega(S^n))\cong\mathbb{Z}[y]$ as well as $H^*(\Omega(S^n))\cong\Gamma[\gamma]$ where $|y|=|\gamma|=n-1$ and $|x|=n$.\\
The term $\mathbb{E}^2$ of the loop homology algebra spectral sequence, and term $E_2$ of the Serre cohomology spectral sequence of the fibration $\Omega(S^n)\hookrightarrow\Lambda(S^n)\rightarrow S^n$ are pictured as follows.\\
\xymatrix@R=2pt@C=5pt{
   &&&&&&&&&&&&&&&&\\
&{xy^3}&&{y^3}&&3(n-1)&&&&&&&{\gamma _3}\ar[1,2]&&{x\gamma _3}&&\\
&{xy^2}&&{y^2}\ar[-1,-2]&&2(n-1)&&&&&&&{\gamma _2}\ar[1,2]&&{x\gamma _2}&&\\
&{xy}&&{y}\ar[-1,-2]&&n-1&&&&&&&{\gamma _1}\ar[1,2]&&{x\gamma _1}&&\\
\ar[0,5] &{x}&&{\bullet}\ar[-1,-2]&&&&&&&&\ar[0,5]&{\bullet}&&{x}&&\\
&-n&&0\ar[-5,0]&\restore&&&&&&&&0\ar[-5,0]&&n&&\\
&&&&&{}\save[]+<0cm,0cm>*\txt<40pc>{\textit{DIAGRAM 1}:\\ Term $\mathbb{E}^2$ of the loop homology algebra spectral sequence and Term $E_2$ of \\the Serre cohomology spectral sequence for the fibration $\Omega(S^n)\hookrightarrow\Lambda(S^n)\rightarrow S^n$}\\
}\\[0.3 cm]
For placement reasons the only nontrivial differentials in either spectral sequence can occur in the $n$-th page. The method of Brown-Shih as described in \cite{McCleary} entails in this case that $S^n, n$ odd is a mod odd $H$-space and therefore \begin{eqnarray*}\Omega(S^n)\hookrightarrow\Lambda(S^n)\rightarrow S^n\end{eqnarray*} is a principal fibration with cross section as described above and therefore the corresponding Serre spectral sequence has to collapse and therefore all differentials in the loop homology algebra sequence have to be zero as well. Our own argument is an elementary computation with the Serre spectral sequence. Since the fibration $\Omega(S^n)\hookrightarrow\Lambda(S^n)\rightarrow S^n$ has a cross section the elements in the base of the corresponding Serre spectral sequence are infinite cycles. The picture in the right shows $d^n(\gamma _1)=0$. Therefore $d_n(y)=0$ has to hold as well and since the differentials fulfill the Leibniz rule and $H_*(\Omega(S^n))$ is a polynomial algebra in one generator the spectral sequence $\{\mathbb{E}^r,d_r \}$ collapses.\\
Note that for the loop homology algebra spectral sequence we always have to deal with two extension issues. The first question is whether $\mathbb{E}^\infty$ is isomorphic as graded $\mathbb{Z}$-module to $\mathbb{H}_*(\Lambda (M))$ and the second extension issue is whether this is an isomorphism of algebras. In this case there are no extension issues because $\mathbb{E}^{\infty}$ is free commutative as abelian group and for the multiplication $\mathbb{E}^{\infty}$ is free commutative and the target commutative. It follows from Theorem 1.6(2) that the multiplicative structure we just computed is the Chas-Sullivan product. We have therefore proved the Theorem 1.2.
\section{Complex projective spaces}
Before we treat the case of even dimensional spheres let us look at the complex projective spaces $\mathbb{C}P^n$.
%Remember from [4] that $\Omega(\mathbb{C}P^n)\simeq S^1\times \Omega(S^{2n+1})$ and for the cohomology Serre spectral sequence we had $H^*(\mathbb{C}P^k)\cong $FIXME%\frac{\mathbb{Z}[x,y]\otimes E[z]}{(x^{k+1}, )}$\\
%The existence of a cross section given by the constant loops implies as above $d^2\equiv 0$ and Leibnitz rules we obtain further $d_2\simeq 0$ for $\mathbb{E}_2$.
Recall that $\Omega(\mathbb{C}P^n)\simeq S^1\times \Omega(S^{2n+1})$ since we can loop the principal fibration $
S^1\rightarrow S^{2n+1}\rightarrow \mathbb{C}P^n$ with classifying map $\mathbb{C}P^n\rightarrow\mathbb{C}P^{\infty}$ and obtain a principal fibration (up to homotopy) $\Omega (S^{2n+1})\rightarrow \Omega (\mathbb{C}P^n)\rightarrow S^1$. Since $\pi _1(\Omega (\mathbb{C}P^n))\cong\mathbb{Z}\cong\pi _1(\mathbb{C}P^{\infty})$ the fibration has a cross section and then \begin{eqnarray*}\Omega(\mathbb{C}P^n)\simeq S^1\times \Omega(S^{2n+1})\end{eqnarray*} as $H-$spaces follows because a principal fibration with cross section is trivial. We therefore have an isomorphism of Hopf algebras \begin{eqnarray*}H_*(\Omega(\mathbb{C}P^n))\cong E(z)\otimes \mathbb{Z}[y]
\end{eqnarray*} with deg$(z)=1$, deg$(y)=2n$. The cohomology spectral sequence for the fibration \begin{eqnarray*}\Omega (\mathbb{C}P^n)\rightarrow\Lambda (\mathbb{C}P^n)\overset{eval}{\longrightarrow }\mathbb{C}P^n\end{eqnarray*} %, where $eval:\Lambda (\mathbb{C}P^n)\rightarrow\mathbb{C}P^n$ denotes evaluation at the basepoint of $S^1$ 
 was completely computed in \cite{[7]} so we simply quote from there. The fibration \begin{eqnarray*}\Omega (\mathbb{C}P^n)\rightarrow\Lambda (\mathbb{C}P^n)\overset{eval}{\longrightarrow }\mathbb{C}P^n\end{eqnarray*} has a cross section. Therefore all elements of the base are infinite cycles, so we know that for $n>1$ the differential $d^2\equiv 0$. For placement reasons the only nonzero differential can occur in the term $E^{2n}$ which is determined by $d^{2n}(y)=z\otimes (n+1)x^n$ as is seen from the following map of fibrations:\\
\xymatrix{ 
{\Omega(\mathbb{C}P^{n})}\ar[1,0]^{}\ar[0,1]&{\Lambda (\mathbb{C}P^{n})}\ar[1,0]^{exp}\ar[0,1]^{eval}&{\mathbb{C}P^{n}}\ar[1,0]^{\Delta}\\
{\Omega(\mathbb{C}P^{n})}\ar[0,1]\restore&{Map(I,\mathbb{C}P^{n})}\ar[0,1]^{eval}&{\mathbb{C}P^{n}\times \mathbb{C}P^{n}},\\
&&&&{}\save[]+<-6cm,0cm>*\txt<20pc>{\textit{DIAGRAM 2}: A Map of Fibrations}\\
} \\[0.3cm]
where \begin{eqnarray*}exp:\Lambda (\mathbb{C}P^n)\rightarrow Map(I,\mathbb{C}P^n)\end{eqnarray*} is the map between function spaces induced by the exponential map $I\rightarrow S^1$ and \begin{eqnarray*}eval:Map(I,\mathbb{C}P^n)\rightarrow\mathbb{C}P^n\times\mathbb{C}P^n\end{eqnarray*} denotes evaluation at the endpoints of $I$.\\
 The preceding discussion allows us to draw the following picture for the term $\mathbb{E}^2$ of the loop homology spectral sequence.\\
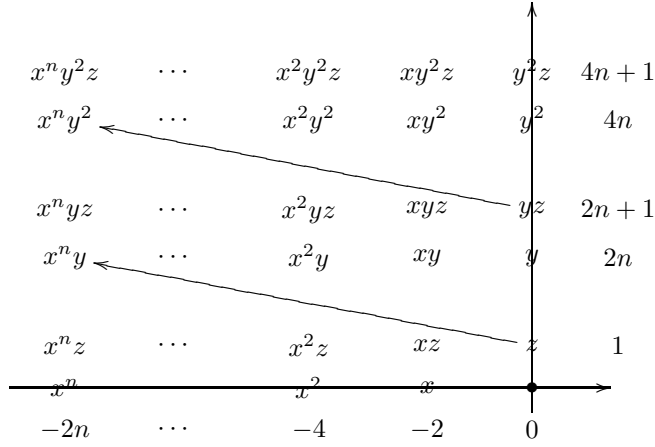
\xymatrix@R=2pt@C=5pt{
&&&&&&&&&&&\\
&&&&&&&&&&&\\
&&&&&&&&&&&\\
&{x^ny^2z}&&{\cdots}&&&{x^2y^2z}&&{xy^2z}&&{y^2z}&4n+1\\
&{x^ny^2}&&{\cdots}&&&{x^2y^2}&&{xy^2}&&{y^2}&4n\\
&&&&&&&&&&&\\
&&&&&&&&&&&\\
&{x^nyz}&&{\cdots}&&&{x^2yz}&&{xyz}&&yz\ar[-3,-9]&2n+1\\
&{x^ny}&&{\cdots}&&&{x^2y}&&{xy}&&y&2n\\
&&&&&&&&&&&\\
&&&&&&&&&&&\\
&{x^nz}&&{\cdots}&&&{x^2z}&&xz&&z\ar[-3,-9]&1\\
\ar[0,11]&{x^n}&&&&&{x^2}&&x&&{\bullet}&\\
&{-2n}&&{\cdots}&&&\restore-4&&-2&&0\ar[-13,0]&\\
&&&&{}\save[]+<2cm,0cm>*\txt<20pc>{\textit{DIAGRAM 3}: Term $\mathbb{E}^{2n}$}\\
&&&&\\
&&&&\\
}\\[0.3cm]From the $\mathbb{E}_2$-page it is not obvious that $d_2(y)=0$, but, again the knowledge of the cohomolgy Serre spectral sequence implies $d_2(y)=0$.
Quoting \cite{[7]} again we obtain with the Universal Coefficient Theorem that $d_{2n}(z)=(n+1)x^n\otimes y$.\\
For placement reasons we have $\mathbb{E}^{2n+1}=\mathbb{E}^{\infty}$.
Knowing all differentials on all generators we obtain 
\begin{eqnarray*}
\underset{p+q=i}{\bigoplus }\mathbb{E}_{p,q}^{\infty}\cong\begin{cases}\mathbb{Z}\text{ if $i=-2n,-2n+1,-2n+2,\cdots ,i\neq 2mn, m\geq 0$}
\\\mathbb{Z}\oplus\mathbb{Z}_{n+1}\text{ if $i=2nm$, $m\geq 0$}
\end{cases}.
\end{eqnarray*}
This is isomorphic to $\mathbb{H}_i(\Lambda (\mathbb{C}P^n);\mathbb{Z})$ for all $i$ as a $\mathbb{Z}-$module
since the module structure does not leave any room for extension issues. For the algebra structure we refer to \cite{[2]}. Let $w=x\otimes z$. The infinite cycles $x\otimes 1$ and $1\otimes y$ represent classes in $\mathbb{H}_{-2}(\Lambda (\mathbb{C}P^n));\mathbb{Z})$ and $\mathbb{H}_{2n}(\Lambda (\mathbb{C}P^n));\mathbb{Z})$ respectively and the ideal $(x^n\otimes 1, (n+1)x^n \otimes y)$ vanishes in $\mathbb{E}_{*,*}^{\infty}$. Therefore we obtain a subalgebra  
\begin{eqnarray*}\frac{ 
 \mathbb{Z}[x\otimes 1,1\otimes y]}{(x^n\otimes 1, (n+1)x^n \otimes y)}\subset \mathbb{E}_{*,*}^{\infty} .\end{eqnarray*}
Note that $w$ is an infinite cycle and represents a class in $\mathbb{H}_{-1}(\Lambda (\mathbb{C}P^n));\mathbb{Z})$ and moreover $w^2=0=w\cdot x^n\otimes 1$. Therefore $\mathbb{E}_{*,*}^{\infty}$ can be written as
\begin{eqnarray*}\mathbb{E}_{*,*}^{\infty}\cong\frac{ E(w)\otimes
 \mathbb{Z}[x\otimes 1, 1\otimes y]}{( x^{n+1}\otimes 1, (n+1)x^n \otimes y,w\cdot x^n)}.\end{eqnarray*}
Careful inspection of the dimensions leaves no possible ambiguities unless $n=1$. Since we have $\mathbb{C}P^1\cong S^2$ we postpone this to the next section. It follows from Theorem 1.6.2 that the multiplicative structure we just computed is the Chas-Sullivan product.
We have thus proved Theorem 1.3.
\section{Even dimensional spheres}
Before we come to our own argument we briefly sketch the method of Brown and Shih as described in \cite{McCleary} which is the basis for the computation of Cohen Jones and Yan in the case of even dimensional spheres. 
Given the fibration \begin{eqnarray*}eval:\Lambda (S^{2k})\rightarrow S^{2k}\end{eqnarray*} with fibre $\Omega S^{2k}$, let $\Omega _{eval}$ denote the pullback over the evaluation map (at $0$), \begin{eqnarray*}eval_0:(S^{2k})^{I}\rightarrow S^{2k}\end{eqnarray*} of the fibration. Let \begin{eqnarray*}\widetilde{eval}:\Lambda (S^{2k})^{I}\rightarrow\Omega _{eval}\end{eqnarray*} be the mapping $\widetilde{eval}(w)=(eval\circ w, w(0))$. The mapping $eval$ is a fibration since there is a section $\sigma :\Omega _{eval}\rightarrow \Lambda (S^{2k})^{I}$ to $\widetilde{eval}$. Notice that $\Omega S^{2k}\times \Omega S^{2k}\subseteq\Omega _{eval}$ which is the subset of $(S^{2k})^I\times \Lambda (S^{2k})$ given by $\{\{(w,e) |w(0)=eval(e)\}$. Restricting $\sigma $ to $\Omega S^{2k}\times \Omega S^{2k}$ we get a mapping to $\Lambda (S^{2k})^I$ such that \begin{eqnarray*}eval\circ \sigma (0)\in S^{2k} \text{ and } eval\circ \sigma (1)\in \Omega S^{2k}.\end{eqnarray*} Composition with evaluation at $1$ this describes an action $\alpha = eval_1 \circ \sigma :\Omega S^{2k}\times \Omega S^{2k}\rightarrow \Omega S^{2k}$ given by $\alpha(w,\lambda)=w^{-1}\lambda w.$
The single differential in the homology Serre spectral sequence for the fibration $\Omega(S^{2k})\hookrightarrow \Lambda(S^{2k})\rightarrow S^{2k}$ can be described with the method of Brown and Shih as
\begin{eqnarray*}
d_{2k}(u\otimes b^j)&=&\alpha_*(\tau(u)\otimes b^j)\\
&=&(-1)^{j|\tau(u)|}b^j\cdot\tau(u)+\chi(\tau(u))\cdot b^j\\
&=&(-1)^{j|\tau(u)|}b^j\cdot\tau(u)-\tau(u)\cdot b^j\\
&\cong &\begin{cases}-2b^{j+1} \text{ for $j$ odd}
\\ 0\text{ for $j$ even}\end{cases}\\
\end{eqnarray*}
where $\tau :H_{2k}(S^{2k})\rightarrow H_{2k-1}(\Omega S^{2k})$ denotes the transgression isomorphism, $\chi$ denotes the canonical antiautomorphism of the Hopf algebra $H_{*}(\Omega S^{2k})$ and $\cdot$ denotes the Pontryagin product on the loop space. Pictorially we obtain\\
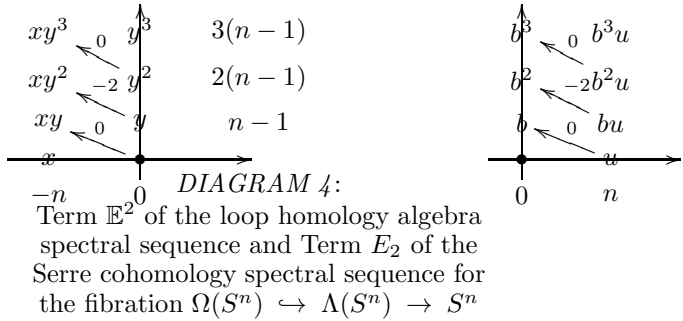
\xymatrix@R=2pt@C=5pt{
   &&&&&&&&&&&&&&&&\\
&{xy^3}&&{y^3}&&3(n-1)&&&&&&&{b^3}&&{b^3u}&&\\
&{xy^2}&&{y^2}\ar[-1,-2]_{0}&&2(n-1)&&&&&&&{b^2}&&{b^2u}\ar[-1,-2]_{0}&&\\
&{xy}&&{y}\ar[-1,-2]_{-2}&&n-1&&&&&&&{b}&&{bu}\ar[-1,-2]_{-2}&&\\
\ar[0,5] &{x}&&{\bullet}\ar[-1,-2]_{0}&&&&&&&&\ar[0,5]&{\bullet}&&{u}\ar[-1,-2]_{0}&&\\
&-n&&0\ar[-5,0]&&&&&&&&& {0}\ar[-5,0]&&n&&\\
&&&&&\restore&&&&&&&&&&&\\
&&&&&{}\save[]+<0cm,0cm>*\txt<20pc>{\textit{DIAGRAM 4}: \\ Term $\mathbb{E}^2$ of the loop homology algebra spectral sequence and Term $E_2$ of the Serre cohomology spectral sequence for the fibration $\Omega(S^n)\hookrightarrow\Lambda(S^n)\rightarrow S^n$}\\
}\\[0.3cm]
where the right hand side indicates the differentials computed via Brown-Shih and the left hand side indicates how Cohen, Jones, and Yan used this information to determine the differentials in the loop homology spectral sequence.\\
Our own method goes as follows.
For even dimensional spheres recall that completely analogously to \cite{[7]} there is map of fibrations:\\
\xymatrix{ 
{\Omega(S^{n})}\ar[1,0]^{}\ar[0,1]&{\Lambda (S^{n})}\ar[1,0]^{exp}\ar[0,1]^{eval}&{S^{n}}\ar[1,0]^{\Delta}\\
{\Omega(S^{n})}\ar[0,1]\restore&{Map(I,S^{n})}\ar[0,1]^{eval}&{S^{n}\times S^{n}}.\\
&&&&{}\save[]+<-6cm,0cm>*\txt<20pc>{\textit{DIAGRAM 5}: A Map of Fibrations}\\
} \\[0.3cm]
%where $eval: \Lambda (S^n)\rightarrow S^n$ denotes evaluation at the base point of $S^1$ and $eval:Map(I,S^n)\rightarrow S^n\times S^n$ denotes evaluation at the two endpoints of the unit intervall $I$.
 Recall that for $n$ even $H_*(\Omega S^n)\cong\mathbb{Z}[y]$, $H^*(\Omega S^n)\cong\Gamma[\gamma]\otimes E(z)$ and $H^*(S^n)\cong E(x)$ with $|y|=n-1$, $|\gamma|=2n-2$, $|z|=n-1$,  and $|x|=n$, 
and that moreover the inclusion of $S^n\hookrightarrow Map(I,S^n)$ as the constant map is a deformation retract. This map of fibrations pictured in DIAGRAM $5$ induces a map of the respective Serre spectral sequences in cohomology. 
For the term $E_2$ of the cohomology spectral sequence of the fibrations \\$\Omega(S^{n})\hookrightarrow S^{n}\rightarrow S^{n}\times S^{n}$ and $\Omega(S^{n})\hookrightarrow \Lambda(S^{n})\rightarrow S^{n}$ we obtain respectively the following diagram.\\
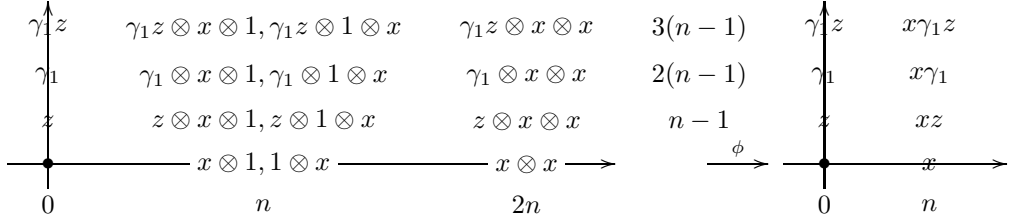
\xymatrix@R=2pt@C=5pt{
   &&&&&&&&&&&\\
&{\gamma _1 z}&&{\gamma _1 z\otimes x\otimes 1,\gamma _1 z\otimes 1\otimes x}&&{\gamma _1 z\otimes x\otimes x}&&3(n-1)&&{\gamma _1 z}&&{x\gamma _1 z}&&\\
&{\gamma _1}&&{\gamma _1\otimes x\otimes 1,\gamma _1 \otimes 1\otimes x}&&{\gamma _1\otimes x\otimes x}&&2(n-1)&&{\gamma _1}&&{x\gamma _1}&&\\
&{z}&&{z\otimes x\otimes 1,z\otimes 1\otimes x}&&{z\otimes x\otimes x}&&n-1&&{z}&&{xz}&&\\
\ar@{-}[0,3]&{\bullet}&&{x\otimes 1,1\otimes x}\ar@{-}[0,2]&&{x\otimes x}\ar[0,1]&&\ar[0,1]^{\phi}&\ar[0,5]&{\bullet}&&{x}&&\\
&0\ar[-5,0]&&n&\restore&2n&&&&0\ar[-5,0]&&n&&\\
&&&&&{}\save[]+<0cm,0cm>*\txt<20pc>{\textit{DIAGRAM 6}: Term $E_2$ and Term $E_2$}\\
}\\[0.3 cm]
For the term $\mathbb{E}_2$ of the the loop homology spectral sequence we also have the following diagram.\\
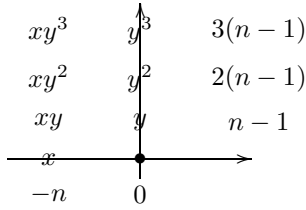
\xymatrix@R=2pt@C=5pt{
   &&&&&&\\
&{xy^3}&&{y^3}&&3(n-1)&\\
&{xy^2}&&{y^2}&&2(n-1)&\\
&{xy}&&{y}&&n-1&\\
\ar[0,5] &{x}&&{\bullet}&&&\\
&-n&&0\ar[-5,0]&\restore&&\\
&&&&&{}\save[]+<0cm,0cm>*\txt<20pc>{\textit{DIAGRAM 7}: Term $\mathbb{E}^2$}\\
}\\[0.3 cm]The next observation is crucial.
Notice that the existence of the Chas Sullivan product implies already that in the spectral sequence $\{\mathbb{E}^r,d_r\}$ not all differentials can be zero because then we had
\begin{eqnarray*}
\mathbb{H}_*(\Lambda(S^{n}))\cong H_*(\Omega (S^{n}))\otimes H^*(S^n)\cong\mathbb{Z}[y]\otimes E(x). 
\end{eqnarray*}
The loop homology is a connected graded commutative algebra thus every subalgebra is graded commutative. But $y\cup y= (-1)\cdot y\cup y$ since the degree of $y$ is odd which would imply $y^2=0$ which is a contradiction. Note that in the case of rational coefficients this argument entirely determines the pattern of the spectral sequence and thus the loop homology algebra since then we have for the preceding argument $d_n(y)=\alpha \cdot xy^2$ for $\alpha\neq 0$ and moreover that because of the graded commutativity $d_n(y^2)=d_n(y^{2k})=0$ for $k\geq 1$ since 
\begin{eqnarray*}d_n(y^{2k})=d_n(y^k\cup y^k)=d_n(y^k)\cup y^k-y^k\cup d_n(y^k)=\alpha\cdot xy^{2k+1}-y^k\cup \alpha\cdot xy^{k+1}=0,
\end{eqnarray*}
and moreover $d_n(y^{2k+1})=\alpha xy^{2k+2}$ for $k\geq 0$ as follows from the Leibniz rule.\\
Let us look at the cohomolgy Serre spectral sequence for the fibration $\Omega(S^{n})\hookrightarrow S^{n}\rightarrow S^{n}\times S^{n}$.\\ Since the target of the spectral sequence does not contain elements of odd degree we obtain $d^n(z)=\pm (x\otimes 1 - 1\otimes x)$. Together with Leibniz rule this gives ker $d^n_{n,n-1}=$span$_{\mathbb{Z}}(x\otimes 1 + 1 \otimes x)$ since we have $d^n(z\otimes (a_0\cdot x\otimes 1 + a_1\cdot 1\otimes x))=(a_0-a_1)x\otimes x$ for $a_0,a_1\in\mathbb{Z}$, (see DIAGRAM 6).\\
This implies $d^n(y)=\pm (x\otimes 1 + 1\otimes x)$. The map of fibrations described in DIAGRAM $5$ induces a map of cohomology Serre spectral sequences which we denote $\phi$. Since we have $d^n(\phi (y))=\phi (d^n(y))=\phi(\pm (x\otimes 1 + 1\otimes x))=\pm2x\otimes z$ the naturality property  tells us $d^n(y)=\pm2\cdot x\otimes z$.\\
\xymatrix@R=2pt@C=5pt{
   &&&&&&&&&&&&&&&&\\
&{xy^3}&&{y^3}&&3(n-1)&&&&&&&{z\gamma _1}\ar[1,2]&&{xz\gamma _1}&&\\
&{xy^2}&&{y^2}\ar[-1,-2]&&2(n-1)&&&&&&&{\gamma _1}\ar[1,2]&&{xz\gamma _1}&&\\
&{xy}&&{y}\ar[-1,-2]&&n-1&&&&&&&{z}\ar[1,2]&&{xz}&&\\
\ar[0,5] &{x}&&{\bullet}\ar[-1,-2]&&&&&&&&\ar[0,5]&{\bullet}&&{x}&&\\
&-n&&0\ar[-5,0]&&&&&&&&& {0}\ar[-5,0]&&n&&\\
&&&&&\restore&&&&&&&&&&&\\
&&&&&{}\save[]+<0cm,0cm>*\txt<20pc>{\textit{DIAGRAM 8}: \\ Term $\mathbb{E}^2$ of the loop homology algebra spectral sequence and Term $E_2$ of the Serre cohomology spectral sequence for the fibration $\Omega(S^n)\hookrightarrow\Lambda(S^n)\rightarrow S^n$}\\
}\\[0.3 cm]
For the loop homology algebra sequence we obtain that the only nonzero differential is $d_{n}(y)=\pm(2\cdot x \otimes y^2)$. Since $\mathbb{E}^{n+1}=\mathbb{E}^{\infty}$ for placement reasons and the differential on all generators are determined 
 we obtain 
\begin{eqnarray*}
\underset{p+q=i}{\bigoplus }\mathbb{E}_{p,q}^{\infty}\cong\begin{cases}\mathbb{Z}\text{ if $i=-n,2m(n-1),-1+2m(n-1), m\geq 0$}
\\\mathbb{Z}_{2}\text{ if $i=-n+2m(n-1), m\geq 1$}
\end{cases}
\end{eqnarray*}
which is isomorphic to $\mathbb{H}_i(\Lambda (S^n);\mathbb{Z})$ for all $i$ as $\mathbb{Z}-$module
since the module structure does not allow any ambiguities. For the multiplicative structure for dimension reasons we can only get an extension problem for $n=2$. In this case we can argument as described in \cite{[2]}. The infinite cycle $1\otimes y^2$ represents a class in $\mathbb{H}_{2(n-1)}(\Lambda (S^n);\mathbb{Z})$. Moreover $1\otimes y^{2j}$ represents a class in $\mathbb{H}_{2j(n-1)}(\Lambda (S^n);\mathbb{Z})$. The ideal $(2x\otimes y^2)$ vanishes in $\mathbb{E}_{*,*}^{\infty}$. Thus $\mathbb{Z}[x\otimes 1,1\otimes y^2]/(x^2\otimes 1,2x\otimes y^2 )$ is a subalgebra of $\mathbb{E}_{*,*}^{\infty}$. Let $x\otimes y\in\mathbb{E}_{-n,n-1}^{\infty} $. We obtain that $\mathbb{E}_{*,*}^{\infty}$ is generated by $\mathbb{Z}[x\otimes 1,1\otimes y^2]/(x^2\otimes 1,2x\otimes y^2 )$ and $x\otimes y$. For dimension reasons we obtain $x^2\otimes 1=x^2\otimes y=0\in \mathbb{E}_{*,*}^{\infty}$, and so
\begin{eqnarray*}
\mathbb{E}_{*,*}^{\infty}\cong \frac{E(x\otimes y)\otimes\mathbb{Z}[x\otimes 1,1\otimes y^2]}{(x^2\otimes 1,x^2\otimes y, 2x\otimes y^2 )}
\end{eqnarray*}
For $n>2$ there cannot be any extension issues for dimensional reasons. It follows from Theorem 1.6.2 that the multiplicative structure we just computed is the Chas-Sullivan product, so
\begin{eqnarray*}
\mathbb{H}_*(\Lambda (S^n);\mathbb{Z}) \cong \frac{E(x\otimes y)\otimes\mathbb{Z}[x\otimes 1,1\otimes y^2]}{(x^2\otimes 1,x^2\otimes y, 2x\otimes y^2 )} \text{ for $n$ even and $n>2$}.
\end{eqnarray*}
.\\
For $n=2$ there is a potential extension problem. The ambiguity of choice of class represented by $1\otimes y^2\in\mathbb{E}_{0,2}^{\infty}$ lies in $\mathbb{E}_{-2,4}^{\infty}$ generated by $x\otimes y^2$. Since we have $x^2\otimes 1=0\in\mathbb{H}_{2}(\Lambda (S^2);\mathbb{Z})$ it follows that any choice $v\in\mathbb{H}_{2}(\Lambda (S^2);\mathbb{Z})$ will satisfy $2x\otimes 1\cdot v=0$. Thus any choice of $v$ together with $x\otimes 1$ and $x\otimes y$ will generate the same algebra, 
\begin{eqnarray*}\mathbb{H}_*(\Lambda (S^2);\mathbb{Z}) \cong\frac{E(x\otimes y)\otimes\mathbb{Z}[x\otimes 1,v]}{(x^2\otimes 1,x^2\otimes y, 2x\otimes 1\cdot v )}.\end{eqnarray*}
 It follows from Theorem 1.6.(2) that the multiplicative structure we just computed is the Chas-Sullivan product.
This proves our last Theorem 1.4.

Dr Nora Seeliger\\
Universität Regensburg, Universitätsstraße 31, 93040 Regensburg, Germany.\\
email: seeliger@math.univ-paris13.fr
\end{document}